# Badly approximable numbers and the growth rate of the inclusion length of an almost periodic function

Anikushin Mikhail
*demolishka@gmail.com*
*Scientific supervisor: Prof. Dr. V. Reitmann, Department of Applied Cybernetics, Faculty of Mathematics and Mechanics, Saint-Petersburg State University*

## Introduction

We study the growth rate of the inclusion length of an almost periodic function. For a given a.p. function such growth rate depends on the arithmetic structure of its Fourier exponents, i.e. on how good they can be approximated by rational numbers. In addition, as appears from the definition, the inclusion length carries some information about the translation numbers (almost periods).

Our result is a lower bound of the growth rate of the inclusion interval of a quasiperiodic function (theorem 3). Here we use methods from dimension theory [2]. We do not assume anything about exponents, but rationally independence. This suggest an idea that this lower bound can be reached (in asymptotic sense) for some "bad" exponents.

Koichiro Naito in his papers on estimates of the fractal dimension of almost periodic attractors [3-4] proved an upper bound of the inclusion length for some class of a.p. functions, using simultaneous Diophantine approximations. For the special case of "badly approximable" exponents, we can see that both estimates (if we consider them as asymptotic estimates) coincide (see theorem 4).

We hope that the ideas and results presented in this paper can be useful not only to understand the nature of badly approximable numbers and almost periods, but also for a more detailed understanding of the structure of almost periodic attractors.

## Main definitions

A subset $A \subset \mathbb{R}$ is called relatively dense in $\mathbb{R}$ if there exists a real number $L > 0$ such that the set $[a, a + L] \cap A$ is not empty for each $a \in \mathbb{R}$. The number L is called an inclusion length.

Let $C_b(\mathbb{R}; \mathbb{C})$ be the space of all bounded continuous functions from $\mathbb{R}$ to $\mathbb{C}$ endowed with the norm

$$||f||_\infty = \sup_{t \in \mathbb{R}} |f(t)|.$$

Consider a function $f \in C_b(\mathbb{R}; \mathbb{C})$. For a given $s \in \mathbb{R}$ define by $f_s$ the $s$-translate of $f$, i.e.

$$f_s(t) = f(t + s) \quad \forall t \in \mathbb{R}.$$

A real number $\tau$ is an $\varepsilon$-almost period for $f$ if $||f_\tau - f||_\infty < \varepsilon$. Denote by $\Omega_\varepsilon(f)$ the set of all $\varepsilon$-almost periods of $f$. The function $f$ is almost periodic if for every

$\varepsilon > 0$ the set $\Omega_\varepsilon(f)$ is relatively dense in $\mathbb{R}$. We denote by $L(\varepsilon)$ the inclusion length of $\Omega_\varepsilon(f)$.

Let $(X, \rho)$ be a compact metric space. Denote by $N_\varepsilon(X)$ the minimum number of balls of radius $\varepsilon$ with centres in $X$ required to cover $X$. The lower and upper limits

$$\underline{\dim}_F(X) = \liminf_{\varepsilon \to 0+} \frac{\ln N_\varepsilon(X)}{\ln(1/\varepsilon)},$$

$$\overline{\dim}_F(X) = \limsup_{\varepsilon \to 0+} \frac{\ln N_\varepsilon(X)}{\ln(1/\varepsilon)}$$

are called the lower fractal dimension and the upper fractal dimension of $X$ respectively [2]. If both values coincide, we use the symbol $\dim_F(X)$ to denote this common value and call it simply the fractal dimension of $X$. We note that the fractal dimension depends on the metric, i.e. for the two topologically equivalent metrics $\rho_1$ and $\rho_2$ on $X$ the limits above can be different. If we want to emphasize the choice of the metric we write $\underline{\dim}_F(X, \rho)$ and $\overline{\dim}_F(X, \rho)$.

Two metrics $\rho_1$ and $\rho_2$ are strongly equivalent if there exists two constants $C_1 > 0$ and $C_2 > 0$ such that
$$C_1 \rho_1(x, y) \leq \rho_2(x, y) \leq C_2 \rho_1(x, y) \quad \forall x, y \in X.$$
It is easy to see that the fractal dimension will not change if we replace the metric by a strongly equivalent one.

### Fractal dimension of the hull of an almost periodic function

Now consider an almost periodic function $f$. The hull of $f$ is defined by the set
$$H(f) = Cl\{f_s \in C_b(\mathbb{R}; \mathbb{C}) \mid s \in \mathbb{R}\},$$
where the closure is taken in the topology of $C_b(\mathbb{R}; \mathbb{C})$. From the Bochner theorem [1] it follows that the hull of an almost periodic function is a compact subset of $C_b(\mathbb{R}; \mathbb{C})$.

For given numbers $a, b \in \mathbb{R}, a < b$, define the set $[a, b]_f = \{f_s \in H(f) \mid s \in [a, b]\}$.

For a given $\varepsilon > 0$ let $L(\varepsilon)$ be an inclusion length of $\Omega_\varepsilon(f)$ and $\delta(\varepsilon)$ be a delta from the definition of the uniform continuity of $f$.

Denote by $N_\varepsilon^{a.p.}$ the minimum number of balls of radius $\varepsilon$ with centres in $\left[-L\left(\frac{\varepsilon}{2}\right), L\left(\frac{\varepsilon}{2}\right)\right]_f$ required to cover $\left[-L\left(\frac{\varepsilon}{2}\right), L\left(\frac{\varepsilon}{2}\right)\right]_f$, i.e. $N_\varepsilon^{a.p.} = N_\varepsilon\left(\left[-L\left(\frac{\varepsilon}{2}\right), L\left(\frac{\varepsilon}{2}\right)\right]_f\right)$. The following lemma holds.

**Lemma 1.**
$$\underline{\dim}_F(H(f)) = \liminf_{\varepsilon \to 0+} \frac{\ln N_\varepsilon^{a.p.}}{\ln(1/\varepsilon)}.$$

$$\overline{\dim}_F(H(f)) = \limsup_{\varepsilon \to 0+} \frac{\ln N_\varepsilon^{a.p.}}{\ln(1/\varepsilon)}.$$

**Idea of the proof.** Note that for each g ∈ H(f) there exists $f_{t_0} \in \left[-L\left(\frac{\varepsilon}{2}\right), L\left(\frac{\varepsilon}{2}\right)\right]_f$ such that $f_{t_0} \in B_\varepsilon(g)$. Every cover by open balls of $H(f)$ corresponds to an open cover by open balls with slightly greater radiuses of $\left[-L\left(\frac{\varepsilon}{2}\right), L\left(\frac{\varepsilon}{2}\right)\right]_f$ of the same cardinality and vice versa. In particular, one can show the inequalities

$$N_{2\varepsilon}^{a.p.} \leq N_\varepsilon(H(f)) \leq N_{\frac{\varepsilon}{2}}^{a.p.}.$$

Thus, the lemma is proved.

Now we can prove an upper estimate of lower and upper fractal dimensions of $H(f)$ in terms of numbers $L(\varepsilon)$ and $\delta(\varepsilon)$.

**Lemma 2.**

$$\underline{\dim}_F(H(f)) \leq \liminf_{\varepsilon \to 0+} \frac{\ln \frac{L(\varepsilon)}{\delta(\varepsilon)}}{\ln(1/\varepsilon)}.$$

$$\overline{\dim}_F(H(f)) \leq \limsup_{\varepsilon \to 0+} \frac{\ln \frac{L(\varepsilon)}{\delta(\varepsilon)}}{\ln(1/\varepsilon)}.$$

**Proof.** Notice that $B_\varepsilon(f_{t_0}) \supset \left[t_0 - \frac{1}{2}\delta\left(\frac{\varepsilon}{2}\right), t_0 + \frac{1}{2}\delta\left(\frac{\varepsilon}{2}\right)\right]_f$. Thus, we can cover $\left[-L\left(\frac{\varepsilon}{2}\right), L\left(\frac{\varepsilon}{2}\right)\right]_f$ by $\frac{2L\left(\frac{\varepsilon}{2}\right)}{\delta\left(\frac{\varepsilon}{2}\right)} + 1$ balls of radius $\varepsilon$. So, $N_\varepsilon^{a.p.} \leq \frac{2L\left(\frac{\varepsilon}{2}\right)}{\delta\left(\frac{\varepsilon}{2}\right)} + 1$. Using lemma 1, we conclude the proof.

### Quasiperiodic case

Within this section, we assume that $f(t) = \sum_{j=1}^n A_j e^{i\lambda_j t}$, where $A_j \in \mathbb{C}$, $A_j \neq 0$, j = 1,2, ..., n and $\lambda_1, ..., \lambda_n \in \mathbb{R}$ are rationally independent (i.e. the equality $p_1\lambda_1 + \cdots + p_n\lambda_n = 0$ for some $p_1, ..., p_n \in \mathbb{Q}$ implies that $p_1 = \cdots = p_n = 0$). Denote by $\mathbb{T}^n$ the n-dimensional plane torus, i.e. $\mathbb{T}^n = \mathbb{R}^n/2\pi\mathbb{Z}^n$ endowed with the metric

$$\rho_{\mathbb{T}^n}(x, y) = \min_{x' \in x, y' \in y} ||x' - y'||.$$

Here $||.||$ is the supremum norm in $\mathbb{R}^n$. Notice that $\rho_{\mathbb{T}^n}(x, y) = \max_{j=1...n} \rho_{\mathbb{T}^1}(x_j, y_j)$, where $x_j$ is the $j$-th coordinate of $x$, i.e. $x = \{(x'_1, ..., x'_n) \in \mathbb{R}^n \mid x'_j \in x_j \; j = 1,2, ..., n\}$. Our aim is to show that $H(f)$ is homeomorphic to $\mathbb{T}^n$ in the strong sense, i.e. the induced metric (w.r.t. the homeomorphism) is strongly equivalent to $\rho_{\mathbb{T}^n}$. At first we need the following theorem [6].

**Theorem 1 (Kronecker's Theorem).** *Suppose that $\lambda_1, \ldots, \lambda_n$ are rationally independent real numbers and $\varkappa_1, \ldots, \varkappa_n$ are arbitrary real numbers; then the system of inequalities*
$$|\lambda_j t - \varkappa_j| < \varepsilon \quad (\mathrm{mod}\ 2\pi) \quad (j = 1, 2, \ldots, n)$$
*has a solution for each $\varepsilon > 0$.*

Now we define the map $\chi: H(f) \to \mathbb{T}^n$ as follows. Consider a function $g \in H(f)$, i.e. $g = \lim_{k \to \infty} f_{t_k}$ for some sequence $t_k \in \mathbb{R}, k = 1, 2, \ldots$. One can find a convergent subsequence $\{t'_k\} \subset \{t_k\}$ such that
$$t'_k \to \varphi_j \left(\mathrm{mod}\ \frac{2\pi}{\lambda_j}\right) \text{ as } k \to \infty \text{ for every } j = 1, 2, \ldots n.$$
Thus, $g(t) = \sum_{j=1}^n A_j e^{i(t\lambda_j + \varkappa_j)}$, where $\varkappa_j \in [0, 2\pi)$. Let's put $\chi(g) = (\varkappa_1, \ldots \varkappa_n)$. It is easy to see that $\chi$ is a continuous map. Using Kronecker's theorem one can show that $\chi$ is a bijection and, due to compactness of $H(f)$, it is a homeomorphism. We took the idea of such a construction from [5].

Now consider the induced metric $\rho'$ on $\mathbb{T}^n$, i.e.
$$\rho'(x, y) = \|\chi^{-1}(x) - \chi^{-1}(y)\|_\infty = \sup_{t \in \mathbb{R}} \left|\sum_{j=1}^n A_j e^{i\lambda_j t}\left(e^{ix_j} - e^{iy_j}\right)\right|,$$
where $(x_1, \ldots, x_n) \in x$ and $(y_1, \ldots, y_n) \in y$ are arbitrary. The following theorem holds.

**Theorem 2.** *Metrics $\rho_{\mathbb{T}^n}$ and $\rho'$ are strongly equivalent.*

**Idea of the proof.** From Kronecker's theorem it follows that $\rho'(x, y) = \sum_{j=1}^n |A_j| |e^{ix_j} - e^{iy_j}|$. We have to show that there exist some constants $C_1 > 0$ and $C_2 > 0$ such that the inequalities $C_1 \rho_{\mathbb{T}^n}(x, y) \leq \rho'(x, y) \leq C_2 \rho_{\mathbb{T}^n}(x, y)$ hold for all $x, y \in \mathbb{T}^n$. The existence of $C_2$ is obvious. To find $C_1 > 0$ we suppose the opposite. Then one can find sequences $x^{(k)}, y^{(k)} \in \mathbb{T}^n, k = 1, 2, \ldots$ such that $\frac{\rho'(x^{(k)}, y^{(k)})}{\rho_{\mathbb{T}^n}(x^{(k)}, y^{(k)})} \leq \frac{1}{k}$ and both, $x^{(k)}$ and $y^{(k)}$, tend to zero as $k$ tends to infinity. For some $(x_1^{(k)}, \ldots, x_n^{(k)}) \in x^{(k)}$ and $(y_1^{(k)}, \ldots, y_n^{(k)}) \in y^{(k)}$ we can see that $\left|e^{ix_j^{(k)}} - e^{iy_j^{(k)}}\right| \geq \frac{1}{2} \rho_{\mathbb{T}^1}(x_j^{(k)}, y_j^{(k)})$ for all sufficiently big numbers $k$. Therefore, we have a contradiction.

It is clear that $\underline{\dim}_F(\mathbb{T}^n, \rho_{\mathbb{T}^n}) = \overline{\dim}_F(\mathbb{T}^n, \rho_{\mathbb{T}^n}) = \dim_F(\mathbb{T}^n, \rho_{\mathbb{T}^n}) = n$. Thus, according to theorem 2, $\dim_F(H(f)) = n$. To get the uniform continuity of $f$ we can put $\delta(\varepsilon) = \frac{\varepsilon}{C}$ for some constant $C > 0$ and thus, using lemma 2, we deduce the following statement.

**Theorem 3.** *Under the above assumptions for the function $f$ we have*
$$L(\varepsilon) \geq \left(\frac{1}{\varepsilon}\right)^{n-1+o(1)},$$

i.e. for every $\gamma > 0$ there exists $\varepsilon_0 > 0$ such that the inequality
$$L(\varepsilon) \geq \left(\frac{1}{\varepsilon}\right)^{n-1-\gamma}$$
holds for all $\varepsilon \in (0, \varepsilon_0)$.

### Relation to badly approximable numbers

A $n$-tuple of real numbers $(\alpha_1, \ldots \alpha_n)$ is called badly approximable if for some constant $C > 0$ and for all positive integers $p$ and $q$ the inequality
$$\max_{j=1\ldots n} \left|\alpha_j - \frac{p}{q}\right| \geq C \left(\frac{1}{q}\right)^{1+\frac{1}{n}}$$
holds. Using theorem 4 from [4] and lemma 5 from [3] we get the following theorem.

**Theorem 4 (K. Naito, [3-4]).** *Consider the function $f(t) = e^{i2\pi t} + \sum_{j=1}^{n} e^{i2\pi\lambda_j t}$, where $\lambda_j \in \mathbb{R}, \lambda_j \neq 0, j = 1, 2, \ldots, n$. Suppose that the inverse exponents $\frac{1}{\lambda_1}, \ldots, \frac{1}{\lambda_n}$ are badly approximable; then for some constant $K > 0$ there exists an inclusion length $L(\varepsilon)$ satisfying the inequality*
$$L(\varepsilon) \leq K \left(\frac{1}{\varepsilon}\right)^n.$$

In my oral talk it was observed that for the case of $f(t) = e^{i2\pi t} + e^{i2\pi\lambda t}$ the value
$$\mathcal{D}i(f) := \limsup_{n \to \infty} \frac{\ln L(\varepsilon)}{\ln 1/\varepsilon}$$
carries an information about the complexity of the trajectory of $f$. In general, the value $\mathcal{D}i(f)$ carries an information about almost periods and the arithmetic structure of the exponents of $f$. It is important that we can estimate this value for the almost periodic solutions of some classes of differential equations where the exponents of the solution are unknown.